\documentclass{article}

\usepackage{amsmath,amssymb,amsthm}
\usepackage{mathtools}
\usepackage{enumitem}
\usepackage{graphicx}
\usepackage{ifpdf}
\ifpdf
\fi

\theoremstyle{thmit} % Numbered and Italic
\newtheorem{thm}{Theorem}[section]
\newtheorem{lem}[thm]{Lemma}
\newtheorem{cor}[thm]{Corollary}
\newtheorem{prop}[thm]{Proposition}
\newtheorem{defn}[thm]{Definition} % not standard

\theoremstyle{thmrm} % Numbered and Roman
\newtheorem{exa}{Example}
\newtheorem*{rem}{Remark}
\newtheorem*{oldproof}{Proof}
\renewenvironment{proof}[1][{}]{\begin{oldproof}[#1]}{\qed\end{oldproof}}

\title{Subrings of $\mathbb{C}$ Generated by Angles}
\author{Jackson Bahr, Arielle Roth}

\providecommand{\abs}[1]{\left\lvert#1\right\rvert}
\providecommand{\card}[1]{\abs{#1}}
\def\ep{\varepsilon}
\def\N{\mathbb{N}}
\def\Z{\mathbb{Z}}
\def\R{\mathbb{R}}
\def\C{\mathbb{C}}
\def\Im{\operatorname{Im}}
\def\Re{\operatorname{Re}}

\begin{document}
 
\maketitle 

\begin{abstract}
    Consider the following inductively defined set.  Given a collection $U$ of unit magnitude complex numbers, and a set initially containing just $0$ and $1$, through each point in the set, draw lines whose angles with the real axis are in $U$.  Add every intersection of such lines to the set.  Upon taking the closure, we obtain $R(U)$.  We investigated for which $U$, $R(U)$ is a ring.

    Our main result holds for when $1 \in U$ and $\card{U} \ge 4$.  If $P$ is the set of real numbers in $R(U)$ generated in the second step of the construction, then $R(U)$ equals the module over $\Z[P]$ generated by the set of points made in the first step of the construction.  This lets us show that whenever the pairwise products of points made in the first step remain inside $R(U)$, it is closed under multiplication, and is thus a ring.
\end{abstract}

\section{Introduction}\label{s:1}

Suppose we are given a collection $U$ of unit length elements of $\C$.  If we have some collection of points in $\C$, we can draw lines through each of them with every angle in $U$ (with respect to the real axis).  In this way we can construct intersections of lines and repeat the process.  Specifically, if we start with $0$ and $1$ in the complex plane and continue this construction forever until it is completed, when is the resulting collection of points a subring of the complex numbers?

Note that even though we are drawing lines, only the intersection points are considered to be constructed.  In \cite{buhler}, Buhler et al.\ motivated this construction with a discussion of origami where two folds can intersect to create a reference point.

\begin{defn}
    Let $p,q,\alpha,\beta \in \C$ with $\abs{\alpha} = \abs{\beta} = 1$.  Define $L_{\alpha}(p)$ to be the line through $p$ with angle $\alpha$.  In other words, $L_{\alpha}(p) \coloneqq p + \R\alpha$.  Define $I_{\alpha, \beta}(p, q) \coloneqq L_{\alpha}(p) \cap L_{\beta}(q)$ when $\alpha \neq \pm \beta$ so that an intersection exists.
\end{defn}

\begin{defn}
    Let $U$ be a set of unit magnitude complex numbers.  Set $S_0 = \{0, 1\}$.  For each $n \in \N$, set $$S_{n+1} = \{I_{\alpha, \beta}(p, q) \mid \alpha, \beta \in U\text{, } p,q \in S_n \text{, and }\alpha \neq \pm\beta\}.$$ We then define $R(U) = \bigcup_{n \in \N} S_n$.
\end{defn}

\begin{defn}
    $T = \{z \in \C \mid \abs{z} = 1\}$ which is viewed as a group under complex multiplication.  $T/\{\pm1\}$ will be used for the collection of angles, since $\alpha$ and $\beta$ are considered equivalent iff $\alpha = \pm \beta$.  Unless otherwise specified, $U \subseteq T/\{\pm 1\}$.
\end{defn}

\begin{defn}
    Given $U \subseteq T/\{\pm 1\}$, we define all elements $z \in R(U)$ of the form $I_{\alpha, \beta}(0,1)$ to be elementary monomials, i.e., length 1 monomials.

    Next, if $m$ is a length $k$ monomial, then $I_{\alpha, \beta}(0, m) \in R(U)$ is a length $k+1$ monomial.  In this way we inductively define monomials.
\end{defn}

\begin{prop}[\cite{buhler}]
    \label{formula}
    We can calculate $I_{\alpha,\beta}(p,q)$ as follows for $p,q \in \C$ and $\alpha \neq \beta \in T/\{\pm 1\}$.
    $$I_{\alpha,\beta}(p,q) = \frac{[\alpha,p]}{[\alpha,\beta]}\beta + \frac{[\beta,q]}{[\beta, \alpha]}\alpha$$
    where $[x,y] = x \bar{y} - y \bar{x}$ and $\bar{z}$ is the complex conjugate of $z$.
\end{prop}

\begin{prop}[\cite{buhler}]
    We list some properties of $I_{\alpha,\beta}(p,q)$ below for $w \in T/\{\pm 1\}$ and $r \in \R$.
    \begin{description}[labelindent=2em] % assuming 2em = \parindent since there's a bug
        \item[(Symmetry)]  $I_{u,v}(p,q) = I_{v,u}(q,p)$
        \item[(Reduction)] $I_{u,v}(p,q) = I_{u,v}(p,0) + I_{v,u}(q,0)$
        \item[(Linearity)] $I_{u,v}(rp+q,0) = rI_{u,v}(p,0) + I_{u,v}(q,0)$
        \item[(Rotation)] For $w \in T/\{\pm 1\}$, $wI_{u,v}(p,q) = I_{wu,wv}(wp,wq)$.
    \end{description}
\end{prop}

\begin{lem}[\cite{buhler}]
    \label{group}
    Let $\card{U} \ge 3$ with $1 \in U$.  Then, $R(U)$ is closed under addition and additive inverses.
\end{lem}

\begin{thm}[\cite{buhler}]
    \label{monomialdecomp}
    Let $\card{U} \ge 3$.  $R(U)$ is the collection of integer linear combinations of monomials.
\end{thm}

\begin{rem}
    Since whenever $\card{U} \ge 3$, $R(U)$ is a group under addition, we need only check closure under multiplication to ensure that $R(U)$ is a ring.  
\end{rem}

The authors of \cite{buhler} then studied the case when $U$ is a group.  Specifically, they took the set of unit magnitude complex numbers $T$ (i.e., the unit circle) and considered it to be a group under complex multiplication.  Then they took the quotient of $T$ by $\{-1, +1\}$.  The result can be viewed as the top half of the unit circle.  By convention, whenever we use $U$, we will refer to $U \subseteq T/\{\pm 1\}$ where the elements are viewed as complex numbers.

\begin{thm}[\cite{buhler}]
    \label{groupring}
    Let $U$ be a subgroup of $T/\{-1, +1\}$ with $\card{U} \ge 3$.  Then, $R(U)$ is a ring.
\end{thm}

In their paper, Buhler et al.\ observed that $R(U)$ maybe be a ring even when $U$ is not a group.  They left the question of what properties $U$ must satisfy exactly for $R(U)$ to be a group open.

\section{Three Angles}\label{s:2}

In order to understand $R(U)$, first we looked at $\card{U} = 3$ with $0 \in \arg(U)$.  We found that $R(U)$ has the structure of a lattice and can be understood in terms of one of the elementary monomials.

\begin{lem}
    Let $U = \{1, u, v\}$.  We claim that $R(U)$ is a lattice in $\C$ with the form $R(U) = \Z + I_{u,v}(0,1) \Z$.
    \begin{proof}
        Set $x = I_{u,v}(0,1)$.  From Lemma \ref{group}, we know that $R(U)$ is a subgroup of $\C$ with addition.  Since $1 \in R(U)$ and $x \in R(U)$, we clearly see that $R(U) \supseteq \Z + x\Z$.
        
        We will prove the other containment with induction.  We know that $S_1 = \{x, 1-x, 0, 1\} \subseteq \Z + x \Z$.  Let $p, q \in S_n$, which is assumed to be in $\Z + x \Z$.  Let $\alpha, \beta \in U$.

        We claim that $z = I_{\alpha,\beta}(p,q) \in \Z + x \Z$.  Since $I_{\alpha,\beta}(p,q) = I_{\alpha,\beta}(p, 0) + I_{\beta, \alpha}(q, 0)$, it suffices to prove that $I_{\alpha,\beta}(a+bx, 0) \in Z + x \Z$.  Further note that 
        \begin{align*}
            I_{\alpha,\beta}(a+bx,0) &= I_{\alpha,\beta}(a,0) + I_{\alpha,\beta}(bx,0) 
            \\
            &= aI_{\alpha,\beta}(1,0) + I_{\alpha,\beta}(bx,0).
        \end{align*}
        by linearity.

        $I_{\alpha, \beta}(1,0) \in S_1$, so $I_{\alpha, \beta}(1,0) = 1$, $0$, $x$, or $1-x$.  There are only four choices since if one of the angles is $0$ radians, the resulting point is $0$ or $1$.  If $\alpha,\beta \neq 1$, then there are two choices left, $\alpha = u, \beta = v$ or $\alpha = v, \beta = u$.  One of these yields the point $x$ and the other yields (by the parallelogram law) $1-x$.  Thus $I_{\alpha, \beta}(a,0) \in \Z + x\Z$.
        
        Next, note that $I_{\alpha,\beta}(bx,0) = bI_{\alpha,\beta}(x,0)$.  Thus it suffices to prove that $I_{\alpha,\beta}(x,0) \in \Z + x\Z$.  We have 6 cases.
        \begin{description}
            \item[($u,v$)] Since $x = ru$ for some $r \in \R$, $I_{u,v}(x,0) = rI_{u,v}(u,0) = 0 \in \Z + x\Z$.
            \item[($v,u$)] $I_{v,u}(x,0)$ is the projection of $x$ on to the line $ru$ in the direction of $v$, but $x \in \R u$, so $I_{v,u}(x,0) = x$.
            \item[($u,1$)] $I_{u,1}(x,0)$ is the projection of $x$ on to the real axis in the direction of $u$.  It is easy to see that this must be $0$, since the line from $0$ (which is on the real axis) extending in the $u$ direction intersects $x$.
            \item[($v,1$)] $I_{v,1}(x,0) = 1$, for a similar reason.  The line extending from $1$ (which is on the real axis) in the $v$ direction intersects $x$.
            \item[($1,u$)] $I_{1,u}(x,0)$ is the line crossing through $x+s$ and $ru$ for $s,r \in \R$, but since $x \in \R u$, this intersection is clearly at $x$.
            \item[($1,v$)] $I_{1,v}(x,0)$ is at $x-1$ which is demonstrated by the fact that $I_{1,v}(x,0) + I_{v,1}(x,0) = x$ and $I_{v,1}(x,0) = 1$.
        \end{description}

        All of these points line in $\Z + x\Z$, so we have shown that $R(U)$ for $\card{U}=3$ is of the form $\Z + x\Z$ where $x = I_{u,v}(0,1)$.
    \end{proof}
\end{lem}

\begin{rem}
    Given $U = \{1, u, v\}$, if we find $u', v'$ such that $I_{u',v'}(0,1) = m + I_{u,v}(0,1)$ for $m \in \Z$ and set $U' = \{1, u', v'\}$, by the above structural result $R(U) = R(U')$.
\end{rem}
Theorem \ref{samelattice} expands on this remark and show exactly when $U$ and $U'$ of size three generate the same lattice.

\begin{thm}
    \label{samelattice}
    Let $I_{u,v}(0,1) = x$ and let $I_{u',v'}(0,1) = y$.  Let $x = a+bi$ and $y = c+di$.  Set $U = \{1,u,v\}$ and $U' = \{1, u',v'\}$.  We claim that $R(U) = R(U')$ if and only if $b = \pm d$ and $a \mp c \in \Z$.
    \begin{proof}
        $\Z + x\Z = \Z + y\Z$ means that $\{m + n x \mid m,n \in\Z\} = \{p + q y \mid p,q \in \Z\}$.  For arbitrarily $m,n \in \Z$, $m + nx \in \{p + qy \mid p,q \in \Z\}$ holds iff $nx \in \Z + y\Z$, which is equivalent to $na+nbi = p + qc + qdi$ for some $p, q \in \Z$.

        In order for this to hold, the imaginary parts must equal: $nbi = qdi$ (for any $n$, there is some $q$).  Thus $d \mid b$ (using $n = 1$).  We can make the same argument swapping $x$ and $y$, which tells us that $b \mid d$, so $b = \pm d$ and thus $n = \pm q$.

        Also, the real parts must be equal: $na - qc = p$ (for any $n$ there are such $p,q$).  Above we determined that $n = \pm q$, so $n(a \mp c) = p$.  Such a $p$ exists for any $n$, so $a \mp c \in \Z$.  We showed that if $\Z + x\Z = \Z + y\Z$, then $b = \pm d$ and $a \mp c \in \Z$.

        Now, if we assume that $b = \pm d$ and $a \mp c \in \Z$, then for any $\Z + x\Z = m + na + nbi$, we have
        \begin{align*}
            m + na + nbi &= m + n(k \pm c) + n(\pm d)i
            \\
            &= (m+nk) \pm nc \pm ndi \in \Z + y\Z.
        \end{align*}
        This shows that $\Z + x\Z \subseteq \Z + y\Z$.  Likewise, $\Z + y\Z \subseteq \Z + x\Z$.

        Since $R(U) = \Z + x\Z$ and $R(U') = \Z + y\Z$, we have that $R(U) = R(U')$ if and only if $b = \pm d$ and $a \mp c \in \Z$, so $\Z + x\Z = \Z + y\Z$.
    \end{proof}
\end{thm}

Now that we understand what form $R(U)$ has for $\card{U} = 3$ with $0 \in \arg(U)$, we can easily show exactly when $R(U)$ is a ring.  The only point that gives any difficulty is $x$, one of the two elementary monomials off of the real line.  If we can square this point and the result lies in $R(U)$, then $R(U) = \Z + x\Z$ must be closed under multiplication.

Now we characterize all $U$ with $0 \in \arg(U)$ and $\card{U} = 3$ such that $R(U)$ is a ring.
\begin{thm}
    \label{quadint}
    Let $U = \{1, u,v\}$ and let $I_{u,v}(0,1) = x$.  $R(U)$ is a ring if and only if $x$ is a (non-real) quadratic integer, i.e., $x$ is the root of some monic integer quadratic polynomial.
    \begin{proof}
        First we will prove that if $x$ is a quadratic integer, then $R(U)$ is a ring.  Note that $R(U) = \Z + x\Z$ where $x = I_{u,v}(0,1)$.  Since $R(U)$ is already a group, we need to show closure under multiplication.  We write $(a+bx)(c+dx) = ac + (bc+ad)x + bdx^2$.  Since $x$ is a quadratic integer, $x^2 = \lambda x + \mu$ for some $\lambda, \mu \in \Z$.  Then,
        \begin{align*}
            (a+bx)(c+dx) &= ac + (bc+ad)x + bd(\lambda x + \mu)
            \\
            &= (ac + bd\mu) + (bc+ad+bd\lambda)x
        \end{align*}
        so in fact $R(U)$ is closed under multiplication.

        Now assume that $R(U)$ is closed under multiplication.  Then $(a+bx)(c+dx) \in \Z + x\Z$, but we can expand this:
        \begin{align*}
            (a+bx)(c+dx) = ac + (bc + ad)x + bd x^2 \in \Z + x\Z
        \end{align*}

        Since $ac + (bc+ad)x \in \Z + x\Z$, we know that $bd x^2 \in \Z + x\Z$ for every $b,d \in \Z$.  In particular, this holds for $b = d = 1$, so $x^2 \in \Z + x\Z$.  In other words, $x$ must be a quadratic integer.  Also, if $x \in \R$, then our $R(U)$ is degenerate, so we need $x \notin \R$.
    \end{proof}
\end{thm}

We can compute the intersection point $x$ in terms of $\arg(u)$ and $\arg(v)$ and rephrase Theorem \ref{quadint}.

\begin{cor}
    Let $\arg(U) = \{0, \theta, \phi\}$ with $\phi < \theta$.  Then $R(U)$ is a ring if and only if $$\frac{\tan \theta}{\tan \theta - \tan \phi} + \frac{\tan \phi \tan \theta}{\tan \theta - \tan \phi}i$$ is a quadratic integer.
    \begin{proof}
        We can see from the following figure that $$(1+w) \tan \phi = h = w \tan \theta$$
        so $w = \frac{\tan \phi}{\tan \theta - \tan \phi}$.
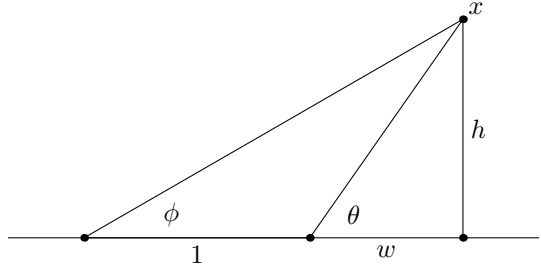
            \begin{center}
                \label{fig1}\includegraphics{figures.1}
            \end{center}
        Immediately, we see also that $h = \frac{\tan \phi \tan \theta}{\tan \theta - \tan \phi}$.  Thus,
        $$x = \frac{\tan \theta}{\tan \theta - \tan \phi} + \frac{\tan \phi \tan \theta}{\tan \theta - \tan \phi}$$
    \end{proof}
\end{cor}

\begin{rem}
    In \cite{nedrenco}, Nedrenco independently characterized $R(U)$ where $\card{U} = 3$, describing $R(U) = \Z + x\Z$ and generalized to when $0 \notin \arg(U)$.  In the same paper, Nedrenco also noted that $R(U)$ is dense when $\card{U} = 4$.  We present what we found independently.
\end{rem}

\section{Four or More Angles}\label{s:3}

Since we understood $R(U)$ for $\card{U} = 3$ in terms of an elementary monomial, we wish to understand $R(U)$ for $\card{U} \ge 4$ in terms of elementary monomials.  Because $R(U)$ is now dense in the complex plane, we cannot hope for an integral basis.  By linearity if we have some $p \in \R \cap R(U)$, then $I_{\alpha, \beta}(0, p) = pI_{\alpha, \beta}(0,1)$.  This means we can scale points.  This motivates our interest in ``projections'' on to the real axis.

\begin{prop}
    Let $U = \{1, u, v, w\}$ with $\arg(u) < \arg(v) < \arg(w) < \pi$.  There are at most eight length-two monomials on the real axis.  There are at most five length-two monomials constructed from elementary monomials of the form $I_{\alpha, \beta}(0,1)$ with $\arg(\alpha) < \arg(\beta)$.  They are $0, 1, x, 1/x, x/(x-1)$ where $x = I_{v, 1}(I_{u, w}(0, 1), 0)$.
    \begin{proof} With the exception of 0 and 1, the only way to construct a length-two monomial on the real axis is to intersect a line through an elementary monomial and the line passing through 0 and 1.  For any given elementary monomial, there are already two lines passing through the point: one passes through 0 and one passes through 1.  Thus there can be at most 6 extra length-two monomials on the real axis, at most three of which created from $z_1, z_2, z_3$ in the form described in the claim, and at most three of which created from $1-z_1, 1-z_2, 1-z_3$ which are of the opposite form.

            
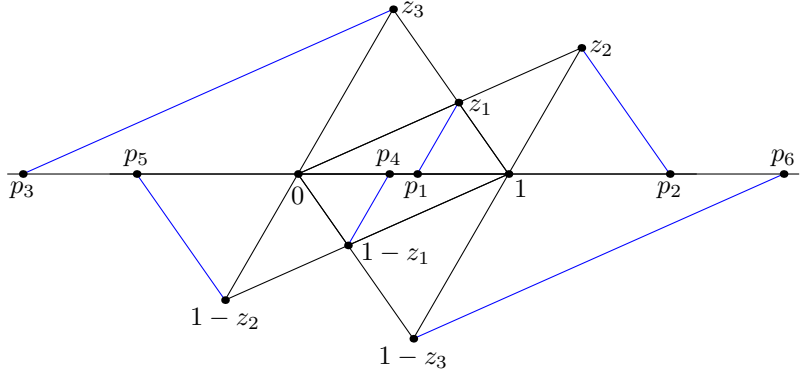
\begin{center}
                \label{fig2}\includegraphics{figures.2}
            \end{center}

            Note that $p_1 = 1 - p_4$, $p_2 = 1 - p_5$, and $p_3 = 1 - p_6$.  As proof, we calculate 
            \begin{align*}
                I_{1, \alpha}(0, I_{\beta, \gamma}(0, 1)) &= I_{1, \alpha}(0, 1 - I_{\gamma, \beta}(0, 1))
                \\
                &= I_{1, \alpha}(0, 1)  - I_{1, \alpha}(0, I_{\gamma, \beta}(0, 1))
                \\
                &= 1 - I_{1, \alpha}(0, I_{\gamma, \beta}(0, 1))
            \end{align*}

            Now we will show that the projections have the described form.  Set $x = p_1$.  Note that the triangle $0-p_1-z_1$ is similar to the triangle $0-1-z_2$, so $\frac{p_1}{1} = \frac{z_1}{z_2}$.  Also, the triangle $0-1-z_1$ is similar to the triangle $0-p_2-z_2$, so $\frac{1}{p_2} = \frac{z_1}{z_2}$.  Thus, $p_2 = 1/x$.
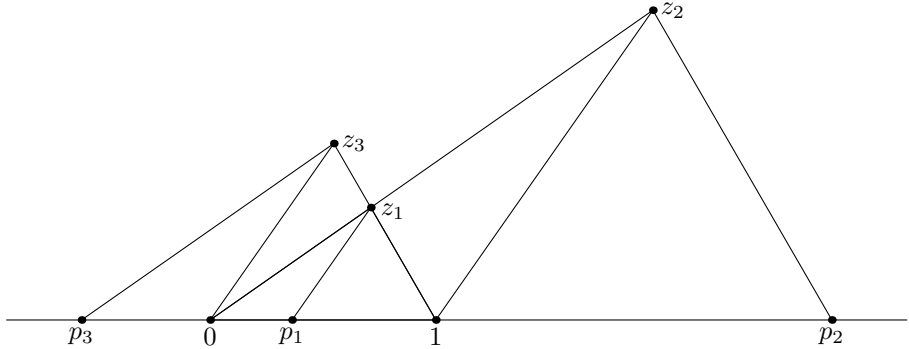
            \begin{center}
                \label{fig3}\includegraphics{figures.3}
            \end{center}

            Next, the triangle $0-p_1 - z_1$ is similar to the triangle $p_3 - 0 - z_3$, so $\frac{\abs{z_1}}{\abs{z_3-p_3}} = \frac{\abs{p_1}}{\abs{p_3}}$.  Also, the triangle $0-1-z_1$ is similar to the triangle $p_3-1-z_3$, so $\frac{\abs{z_1}}{\abs{z_3-p_3}} = \frac{1}{\abs{1-p_3}}$.  We conclude
            \begin{align*}
                \frac{\abs{x}}{\abs{p_3}} &= \frac{1}{\abs{1-p_3}}
                \\
                \abs{p_3-1}\abs{x} &= \abs{p_3}
                \\
                \abs{p_3} &= \abs{\frac{x}{x-1}}
            \end{align*}

            To remove the absolute value signs, we note that since $\arg(z_3) > \arg(z_1)$, the line through $z_3$ with angle $\arg(z_1)$ must intersect the negative real axis, so $p_3 < 0$.  Furthermore, since $x < 1$, $\frac{x}{x-1} < 0$, so we deduce that $p_3 = x/(x-1)$.

    \end{proof}
\end{prop}

Now that we understand a small amount of $\R \cap R(U)$, we can quickly construct an entire ring inside $\R \cap R(U)$ with the scaling mentioned earlier.  Later we will show that what we construct next is exactly $\R \cap R(U)$

\begin{prop}\label{projconstruct}
    Let $0 \in \arg(U)$ with $\card{U} \ge 4$.  Let $P$ be the set of length-two monomials on the real axis.  For any $x \in R(U)$ and any $p \in P$, $px \in R(U)$.  As a result, the ring $\Z[P]$ is constructible, i.e., $\Z[P] \subseteq R(U)$.
    \begin{proof}
        Let $p$ be a projection.  Since $R(U)$ is the collection of finite linear combinations of monomials, it suffices to construct $pm$ for a given monomial $m$, since if we have $x \in R(U)$, we can simply represent $x = \sum_{i=1}^{n} c_i m_i$ for $c_i \in \Z$ and then write $px = \sum_{i=1}^{n} c_i (pm_i)$.  
        
        The proof that $pm \in R(U)$ follows from linearity.  Formally, we rely on induction.
        \begin{description}
            \item[Base Case:] The length of $m$ is one, so $m = I_{\alpha, \beta}(0,1)$ for some $\alpha, \beta \in U$.  Then, $pm = I_{\alpha, \beta}(0,p)$ by linearity, which is in $R(U)$ since $p \in R(U)$.
            \item[Inductive Step:] Suppose every length $n-1$ monomial satisfies the claim.  Let $m$ be of length $n$.  Then, $m = I_{\alpha, \beta}(0, q)$ for some length $n-1$ monomial $q$.  By linearity, $pm = I_{\alpha, \beta}(0, pq)$ which is constructible since $pq \in R(U)$ by the inductive hypothesis.
        \end{description}

        Thus every monomial can be arbitrarily multiplied by projections, so in fact everything in $R(U)$ can be arbitrarily multiplied by projections.  In particular, so can the projections themselves.  This means that arbitrary powers of projections are in $R(U)$.  Furthermore, since $R(U)$ is a group under addition, $\Z[P] \subseteq R(U)$.
    \end{proof}
\end{prop}

\begin{rem}
Since the above result does not rely on the previous two results, this holds even when $\card{U} > 4$.
\end{rem}

Our current goal is to characterize all monomials in terms of $\Z[P]$ and elementary monomials.  By Theorem \ref{monomialdecomp}, if the monomials have a nice enough form, we will be able to understand all of $R(U)$.  Characterizing all monomials starts with the length two monomials.  First, however, we need a quick lemma.

\begin{lem}
    Let $0,\alpha,\beta \in \arg(U)$.  Let $p, q \in R(U)$, and let $x = I_{\alpha, \beta}(p,q)$ and $y = I_{\beta, \alpha}(p,q)$.  Then, $x = p+q-y$.
    \begin{proof}
        Since the lines from $x$ to $q$ and from $p$ to $y$ are parallel, and also the lines from $x$ to $p$ and from $q$ to $y$ are parallel, this forms a parallelogram.  It is clear that $0$, $x-q$, $p-q$, and $y-q$ form a parallelogram and that $x-q + y-q = p-q$, so $x+y-q = p$.  
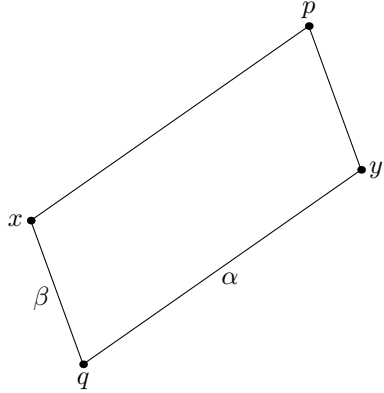
            \begin{center}
                \label{fig4}\includegraphics{figures.4}
            \end{center}
    \end{proof}
\end{lem}

\begin{lem}
    Let $\card{U} \ge 4$ and let $0 \in \arg(U)$.  Let $P$ be the set of projections from the elementary monomials to the real axis along angles in $U$.  Every length two monomial is a $\Z[P]$-linear combination of elementary monomials.
    \begin{proof}
        Let $z = I_{\alpha, \beta}(0,1)$ for some $\alpha,\beta \in U$ and let our length two monomial $m = I_{\gamma, \delta}(0,z)$.  We will prove that $m$ is a $\Z[P]$-linear combination of elementary monomials by cases.

        \begin{description}
            \item[($\delta  = 1$):]
                Note that $$I_{\gamma,0}(0,z) + I_{0,\gamma}(0,z) = z,$$ so $I_{\gamma, \delta}(0,z) = z - I_{0, \gamma}(0,z)$.  Since $I_{0, \gamma}(0,z) \in P$, $m$ is a $\Z[P]$-linear combination of elementary monomials.

            \item[($\delta = \alpha$):]
                Since the line through $z = I_{\alpha, \beta}(0,1)$ with angle $\arg(\alpha)$ passes through the origin, $m = I_{\gamma, \alpha}(0, z) = 0$.  This is trivially a $\Z[P]$-linear combination of elementary monomials.

            \item[($\delta = \beta$):]
                Since the line through $z = I_{\alpha, \beta}(0,1)$ with angle $\arg(\beta)$ passes through $1$, $m = I_{\gamma, \beta}(0,z) = I_{\gamma, \beta}(0,1)$, which is an elementary monomial.

            \item[($\delta \in U \setminus\{1,\alpha,\beta\}$):]
                Let $p = I_{0, \gamma}(0, z)$ be the projection from $z$ to the real axis in the direction of $\gamma$.  Note that $I_{\gamma, \delta}(0, p) = p I_{\gamma, \delta}(0,1)$ by linearity.

                Set $x = I_{\gamma, \delta}(0, p)$.  We that $x + z - p = m$, and since $x = p I_{\gamma, \delta}(0,1)$, this is enough to prove that $m$ is a $\Z[P]$-linear combination of elementary monomials.  Restated, the claim is that
                \begin{align*}
                    I_{\gamma, \delta}(0, I_{0, \gamma}(0,z)) + z - I_{0,\gamma}(0,z) = I_{\gamma, \delta}(0, z)
                \end{align*}
                To prove this, we will show that $I_{\gamma, \delta}(x,z) = m$.  This follows by the fact that $x \in \R \gamma$, so the line through $x$ with angle $\arg{\delta}$ also passes through 0 and thus $I_{\gamma, \delta}(x,z) = I_{\gamma, \delta}(0,z) = m$.

                Furthermore, $I_{\delta, \gamma}(x,z) = p$.  To see this, first note that $I_{\gamma, 0}(z, 0) = p$.  Also, $I_{\delta, 0}(x, 0) = p$, because
                \begin{align*}
                    I_{\delta, 0}(x,0) &= I_{\delta, 0}(I_{\gamma, \delta}(0, p), 0)
                \end{align*}
                and both $x$ and $p$ lie along the same line through $p$ with angle $\arg(\delta)$ (by construction of $x$).

                This means that $x$ and $z$ lie on opposite corners of a parallelogram which has a corner at $p$ through the real axis and another corner through $m$.  Thus, $0$, $(x-p)$, $(z-p)$, and $(m-p)$ form the corners of a parallelogram and $(x-p) + (z-p) = m-p$ so $x+z-p = m$, concluding the proof.
        \end{description}
        Since in all cases $m$ is a $\Z[P]$-linear combination of elementary monomials, we know that every length two monomial is of this form.
    \end{proof}
    
\end{lem}

Now that we understand length two monomials, we can apply induction to characterize all monomials, and thus all of $R(U)$.

\begin{thm}
    Let $0 \in \arg(U)$.  Let $P$ be the set of projections of elementary monomials along lines with angles from $\arg(U)$ on to the real axis.  Then, every monomial in $R(U)$ is a $\Z[P]$-linear combination of elementary monomials.  Indeed, $R(U)$ is the set of $\Z[P]$-linear combinations of elementary monomials.
    \begin{proof}
        We will prove this by induction on the length of the monomial.  Length one monomials are already elementary and length two monomials follow from the above theorem.  Let $m$ be length $n$ and suppose that all length $n-1$ monomials are of this form.  Then,
        \begin{align*}
            m &= I_{\alpha, \beta}(0, m')
            \\
            &= I_{\alpha, \beta}(0, \sum_{i=1}^{k} c_i z_i)
            \\
            &= \sum_{i=1}^{k} c_i I_{\alpha, \beta}(0, z_i)
            \\
            &= \sum_{i=1}^k \left(c_i  \sum_{j=1}^{\ell} d_i x_i\right)
        \end{align*}
        using linearity and the fact that all length two monomials are of this form.  The $c_i$ and $d_i$ are in $\Z[P]$ and the $x_i$ and $z_i$ are elementary monomials.  After simplification, it is easy to see that $m$ is in fact a $\Z[P]$-linear combination of elementary monomials.

        Since everything in $R(U)$ is an integer linear combination of monomials, everything in $R(U)$ is a $\Z[P]$-linear combination of elementary monomials.

        Furthermore, since $\Z[P]$ is constructible by Proposition \ref{projconstruct}, and $p R(U) \subseteq R(U)$ for all $p \in P$, we can construct every $\Z[P]$-linear combination of elementary monomials.  Thus, $R(U)$ equals the set of $\Z[P]$-linear combinations of elementary monomials.
    \end{proof}
\end{thm}

\begin{rem}
    We can alternatively say that $R(U)$ is a $\Z[P]$-module in $\C$ generated by the elementary monomials.
\end{rem}

As in the three-angle case, understanding the structure of $R(U)$ led us to understand when $R(U)$ is a ring in terms of products of elementary monomials.  In fact Theorem \ref{quadint} could probably be seen as a special case of the following theorem.

\begin{thm}
    \label{bigresult}
    Let $U$ with $\card{U} \ge 4$ and $0 \in \arg(U)$ and let $P$ represent the collection of projections.  $R(U)$ is a ring if and only if every pairwise product of elementary monomials is a $\Z[P]$-linear combination of elementary monomials.
    \begin{proof}
        First note that $R(U)$ equals the collection of $\Z[P]$-linear combinations of elementary monomials.  We that the $\Z[P]$-linear combinations of elementary monomials are closed under multiplication if and only if every pairwise product of elementary monomials is a $\Z[P]$-linear combination of elementary monomials.

        Assume that every pairwise product of elementary monomials is as above.  Then, for any $x, y \in R(U)$, we write $x = \sum_{i=1}^{n} c_i x_i$ and $y = \sum_{j=1}^{m} d_j y_j$ for $c_i,d_j \in \Z[P]$ and $x_i,y_j$ elementary monomials.

        Then, $xy = \sum_{i,j} c_i d_j x_i y_j$.  Since $x_i y_j$ is a $\Z[P]$-linear combination of elementary monomials, so is $xy$.  Thus $R(U)$ is a ring.

        Now, suppose that $R(U)$ is a ring.  It must be closed under multiplication, so the pairwise product of elementary monomials must be in $R(U)$, but $R(U)$ is the $\Z[P]$-linear combinations of elementary monomials, so the claim holds.
    \end{proof}
\end{thm}

Since we have at least one projection $p \in (0,1)$, we can construct points close to zero.  Because elements of $R(U)$ scaled by $p$ are still in $R(U)$ and $R(U)$ is a group, it is actually dense in $\C$ as we will prove below.
\begin{thm}
    \label{density}
    If $1 \in U$ and $\card{U} \ge 4$, then $R(U)$ is dense in $\C$.
    \begin{proof}

        Since $R(U)$ is the set of $\Z[P]$-linear combinations of elementary monomials, if $z$ is a non-real elementary monomial and $p \in \Z[P] \cap (0,1)$, we can construct $p^n$ and $p^n z$ which go to zero from, two different directions.

        Let $\ep > 0$ and let $x \in \C$.  Since $R(U)$ is a group under addition, we can construct $a p^{N_1} + b p^{N_2} z$ for all $N_1, N_2 \in \N$.
        
        Since $p \in (0,1)$, we can find $N_2$ such that $\abs{\Im(z) p^{N_2}} < \ep/2$.  To simplify the following expression, write $\theta = \Im(z) p^{N_2}$.  Then there exists a unique $b \in \Z$ such that
        \begin{align*}
            b-1 \le \frac{\Im(x)}{\theta} \le b
        \end{align*}
        So we can show that
        \begin{align*}
            \abs{b \Im(z) p^{N_2}i - \Im(x)i} = \abs{b\theta - \Im(x)} \le \ep/2
        \end{align*}

        Likewise we can find $a, N_1$ such that $\abs{a p^{N_1} - \left(\Re(x) - b p^{N_2} \Re(p)\right)} < \ep/2$.  Once we have such $a \in \Z$ and $N_1 \in \N$, we have
        \begin{align*}
            \abs{a p^{N_1} + b p^{N_2} z - x} &= \abs{a p^{N_1} + b p^{N_2} \Re(z) - \Re(x) + b p^{N_2} \Im(z)i - \Im(x)i}
            \\
            &\le \abs{a p^{N_1} + b p^{N_2} \Re(z) - \Re(x)} + \abs{b p^{N_2} \Im(z) - \Im(x)} 
            \\
            &< \ep
        \end{align*}

        Since $a p^{N_1} + bp^{N_2}z \in R(U)$, and this holds for any $x \in \C$ and for every $\ep > 0$, we can always find a point in $R(U)$ arbitrarily close to any point of $\C$.  Thus, $R(U)$ is dense in $\C$.
        
    \end{proof}
\end{thm}

\section{Some $U$ for Which $R(U)$ Is a Ring}\label{s:4}

Now we can use Theorem \ref{bigresult} to prove that $R(U)$ is a ring for a particular example of $U$.

\begin{exa}
    Let $U = \{1, e^{i\pi/6}, e^{i\pi/3}, e^{i\pi/2}\}$.  $R(U)$ is a ring.
    \begin{proof}
        It suffices to show that all products of elementary monomials are $\Z[P]$-linear combinations of elementary monomials.  Our elementary monomials are $0, 1, z_1, z_2, z_3, 1-z_1, 1-z_2, 1-z_3$, where
    \begin{align*}
        z_1 &= \frac{2\sqrt{3}}{3} e^{i\pi/6}
        \\
        z_2 &= \sqrt{3} e^{i\pi/6}
        \\
        z_3 &= 2 e^{i\pi/3}
    \end{align*}

    First we calculate the projections and get $2/3, 3/2, -2$.  Note that $\Z[2/3, 3/2, -2] = \Z[2/3,3/2] = \Z[1/3,1/2] = \Z[1/6]$.

We calculate all pairwise products of $z_1, z_2, z_3$, since calculating more would be redundant, as the others are either $0$, $1$, or an integer linear combination of $\{1, z_1, z_2, z_3\}$.
    \begin{align*}
        z_1^2 &= \frac{4}{3} e^{i\pi/3} = \frac{2}{3} z_3
        \\
        z_1 z_2 &= 2 e^{i\pi/3} = z_3
        \\
        z_1 z_3 &= \frac{4}{\sqrt{3}} e^{i\pi/2} = \frac{4i}{\sqrt{3}} = 4(z_1-1)
        \\
        z_2^2 &= z_1^2 \frac{z_2^2}{z_1^2} = \frac{9}{4} \cdot \frac{2}{3} z_3 = \frac{3}{2} z_3
        \\
        z_2 z_3 &= \frac{z_2}{z_1} z_1 z_3 = 6(z_1 - 1)
        \\
        z_3^2 &= z_1 z_2 z_3 = 6(z_1^2 - z_1) = 4z_3 - 6z_1
    \end{align*}

    These are all in $R(U)$, so $R(U)$ is closed under multiplication and is a ring.
\end{proof}
\end{exa}

\begin{rem}
    We suspected that perhaps any subset $U$ of a finite group containing a generator for that finite group would result in a ring.  The following example shows that this cannot be necessary.
\end{rem}

\begin{exa}
    Let $U = \{1, e^{i\pi/6}, e^{i\pi/4}, e^{i\pi/3}\}$.  $R(U)$ is a ring.
    \begin{proof}
        As above, it suffices to show that the products of all elementary monomials are $\Z[P]$-linear combinations of elementary monomials.  We go by the convention that $z_1 = I_{e^{i\pi/6}, e^{i\pi/2}}(0, 1)$, $z_2 = I_{e^{i\pi/6}, e^{i\pi/3}}(0, 1)$, and $z_3 = I_{e^{i\pi/3}, e^{i\pi/2}}(0, 1)$ and that $p_1, p_2, p_3$ are projections from $z_1, z_2, z_3$ to the real axis.

        We calculated
        \begin{align*}
            z_1 z_2 &= p_3(1- z_3)
            \\
            z_1 z_3 &= -p_1 z_2 - (p_2 p_3)z_3 + 2 p_3
            \\
            z_2 z_3 &= -p_3 z_2 - (p_2 p_3)z_3 + 2 p_2 p_3
            \\
            z_1^2 &= p_1 p_3(1-z_3)
            \\
            z_2^2 &= p_2 p_3 (1-z_3)
            \\
            z_3^2 &= -6z_2 - 3p_2p_3 z_3 + 3p_3
        \end{align*}
    \end{proof}

\end{exa}

\begin{rem}
    We then suspected that any subset of a finite group might result in a ring.  Our next result shows this too cannot be necessary.
\end{rem}

\begin{exa} 
    Let $U = \{1, e^{i}, e^{2i}, e^{3i}\}$.  $R(U)$ is a ring.
\end{exa}
This example is a special case of Theorem \ref{4theta}.

\begin{rem}
    We strongly suspect that $R(\{1, e^{i\pi/5}, e^{i\pi/4}, e^{i\pi/3}\})$ is not a ring, so we suspect that it is not sufficient for $U$ to just be a subset of a finite group.
\end{rem}

\begin{thm}
    \label{4theta}
    Let $U = \{1, \alpha, \alpha^2, \alpha^3\}$.  $R(U)$ is a ring.
    \begin{proof}
        Set $z_1 = I_{\alpha, \alpha^3}(0,1)$, $z_2 = I_{\alpha, \alpha^2}(0,1)$, and $z_3 = I_{\alpha^2, \alpha^3}(0,1)$.  Since the only elementary monomials are $0,1,z_1, z_2, z_3, 1-z_1, 1-z_2, 1-z_3$, it suffices to check pairwise products of $\{z_1, z_2, z_3\}$.  

Set $p_1 = I_{1, \alpha^2}(0,z_1)$, $p_2 = I_{1, \alpha^3}(0, z_2)$, and $p_3 = I_{1, \alpha}(0, z_3)$.  Then $\Z[P] = \Z[p_1, p_2, p_3]$, since the other projections are 0, 1, $1-p_1$, $1-p_2$, and $1-p_2$.

        First we claim that $z_1 z_2 = z_3$.  We will prove this by calculation.
        \begin{align*}
            z_1 z_2 &= \frac{[1,\alpha^3]}{[\alpha, \alpha^3]} \frac{[1, \alpha^2]}{[\alpha, \alpha^2]} \alpha^2
            \\
            &= \frac{e^{-3i\theta} - e^{3i\theta}}{e^{-2i\theta} - e^{2i\theta}} \frac{e^{-2i\theta} - e^{2i\theta}}{e^{-i\theta} - e^{i\theta}}
            \\
            &= \frac{[1,\alpha^3]}{[\alpha, \alpha^2]}\alpha^2 = \frac{[1,\alpha^3]}{[\alpha^2,\alpha^3]} \alpha^2
            \\
            &= z_3
        \end{align*}

        Next we claim that $z_1/z_2 = p_1$ and $z_2 / z_1 = p_2$.  These can also be calculated but a geometrical figure makes it clear.  
        
        The first claim follows from the fact that the triangles $0 - p_1 - z_1$ and $0 - 1 - z_2$ are similar.  The second claim follows from the similarity of the triangles $0 - 1 - z_1$ and $0 - p_2 - z_2$.

            
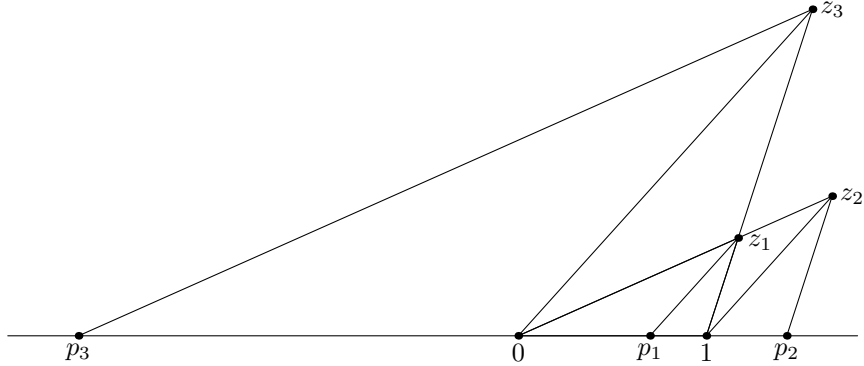
\begin{center}
                \label{fig5}\includegraphics{figures.5}
            \end{center}

        So far we can construct the following pairwise products of elementary monomials.
        \begin{align*}
            z_1^2 &= z_1 z_2 \frac{z_1}{z_2} = p_1 z_3
            \\
            z_1 z_2 &= z_3
            \\
            z_2^2 &= p_2^2 p_1 z_3 = p_2 z_3
        \end{align*}
        We need only construct $z_3^2$ and $z_2 z_3$ since $z_1 z_3 = p_1 z_2 z_3$.

        First we show $z_3^2 = p_3^2(z_3 - z_2)$ algebraically.  We calculated $z_3^2$ using the formula given in Proposition \ref{formula} and obtained
        \begin{align*}
            z_3^2 = 1 + 2 \alpha^2  + 3 \alpha^4 + 2 \alpha^6 + \alpha^8
        \end{align*}
        which is exactly what we found by calculating $p_3^2 (z_3 - z_2)$, so the two must be equal.

        Likewise, we calculated $z_2 z_3$ to be
        \begin{align*}
            z_2 z_3 = 1 + 2 \alpha^2  + 2 \alpha^4 + \alpha^6
        \end{align*}
        which precisely equals $p_3 (1- z_3)$.

Thus all 6 pairwise products of $\{z_1, z_2, z_3\}$ are $\Z[P]$-linear combinations of elementary monomials, so $R(U)$ is a ring.
    \end{proof}
\end{thm}

This characterization of $R(U)$ makes finding examples of rings $R(U)$ a matter of verifying that finitely many products are contained in $R(U)$.  However, finding counterexamples is more difficult.  Some $U$ that are difficult to work with, like $\{1, e^{i}, e^{2i}, e^{3i}\}$, yield rings.  Other $U$ that are nicer to work with, such as $\{1, e^{i\pi/5}, e^{i\pi/4}, e^{i\pi/3}\}$ are suspected to not yield rings.

\section{Open Questions}\label{s:5}

Some open questions we considered in research are posed below.

\begin{enumerate}
    \item How does $1 \notin U$ affect our current results?  Can we still express $R(U)$ as a module over some ring generated by elementary monomials?
    \item When exactly are the products of elementary monomials $\Z[P]$-linear combinations of elementary monomials?
    \item Is $R(\{1, e^{i\pi/5}, e^{i\pi/4}, e^{i\pi/3}\})$ a ring?
    \item What subrings of $\C$ are of the form $R(U)$ for some $U$?
    \item Given $p \in \C$, for which $U$ is $p \in R(U)$?
    \item We can write $I_{u,v}(p,q) = \frac{[u,p]}{[u,v]}v + \frac{[v,q]}{[v,u]}u$ where $[x,y] = x \bar{y} - y \bar{x}$.  Note that $[x,y]$ is an alternating bilinear map.  If $V$ is some vector space equipped with $[\cdot, \cdot]$, an alternating bilinear map into $\R$ and we have some $S \subseteq V$ of allowable ``angles'', we can define $I : S^2 \times V^2 \to V$ via 
\begin{align*}
    I_{u,v}(p,q) = \frac{[u,p]}{[u,v]}v + \frac{[v,q]}{[v,u]}u
\end{align*}
Do similar results hold for this generalization?  Perhaps we could require $V$ to be a normed vector space and say that $S$ is the sphere of radius one.
\end{enumerate}

\bibliographystyle{line}

\end{document}